\documentclass[12pt]{article}

\newcommand{\set}[1]{\{#1\}}

\usepackage{etex}
\usepackage{array}
\usepackage{amsmath, amsthm, amssymb}
\usepackage{pictex}
\usepackage{pstricks}
\usepackage{pst-node}
\usepackage{multirow}
\usepackage{fullpage}
\usepackage{graphicx}

\newtheorem{theorem}{Theorem}[section]

\newtheorem{corollary}[theorem]{Corollary}
\newtheorem{conjecture}{Conjecture}

\newenvironment{definition}[1][Definition]{\begin{trivlist}
\item[\hskip \labelsep {\bfseries #1}]}{\end{trivlist}}
\newenvironment{example}[1][Example]{\begin{trivlist}
\item[\hskip \labelsep {\bfseries #1}]}{\end{trivlist}}

\begin{document}

\renewcommand{\arraystretch}{1}

\title{Spanning Trees in Grid Graphs}   % type title between braces
\author{Paul Raff}         % type author(s) between braces
\date{July 25, 2008}    % type date between braces
\maketitle

\begin{abstract}
A general method is obtained for finding recurrences involving the number of spanning trees of grid graphs, obtained by taking the graph product of an arbitrary graph and path or cycle. The results in this paper extend the work by Desjarlais and Molina and give concrete methods for finding the recurrences. Many new recurrences are found, yielding conjectures on the order of the linear recurrences of grid graphs and graphs obtained by taking the product of a complete graph and a path.
\end{abstract}

\section{Introduction}

The Matrix Tree Theorem of Kirchhoff, a generalization of Cayley's Theorem from complete graphs to arbitrary graphs \cite{Stanley}, gives the number of spanning trees on a labeled graph as a determinant of a specific matrix. If $A = (a_{ij})$ is the adjacency matrix of a graph $G$, then the number of spanning trees can be found by computing any cofactor of the Laplacian matrix of $G$, or specific to the $(n,n)$-cofactor:

\renewcommand{\arraystretch}{2}

\[
\mbox{Number of spanning trees of $G$ } =
\begin{vmatrix}
a_{12} + \ldots + a_{1n} & -a_{12} & \cdots & -a_{1,n-1} \\
-a_{21} & a_{21} + \cdots + a_{2n} & & -a_{2,n-1} \\
\vdots & & \ddots & \vdots \\
-a_{n-1,1} & -a_{n-1,2} & \cdots & a_{n-1,1} + \cdots + a_{n-1,n}
\end{vmatrix} \]

\renewcommand{\arraystretch}{1}

Since determinants are easy to compute, then the Matrix Tree Theorem allows for the computation for the first few numbers in the sequence of spanning trees for families of graphs dependent on one or more parameters. However, the downside of the Matrix Tree Theorem is that it can only produce a sequence of numbers, and cannot \emph{a priori} assist in finding out the recurrence involved with said sequence. In this paper, the motivation is the following families of graphs:
\begin{enumerate}
\item $k \times n$ grid graphs, with $n \to \infty$.
\item $k \times n$ cylinder graphs, with $n \to \infty$.
\item $k \times n$ torus graphs, with $n \to \infty$.
\end{enumerate}

All of the families of graphs mentioned above can be placed into a more general class of graphs of the form $G \times P_n$ or $G \times C_n$, where $P_n$ and $C_n$ denote the path and cylinder graph on $n$ vertices, respectively. For each of these classes, a general method is obtained for finding recurrences for all of the above families of graphs, and explicit recurrences are found for many cases. The only drawback, as it stands, is the amount of computational power needed to obtain these recurrences, as the recurrences are obtained through characteristic polynomials of large matrices. The result is at least 15 new sequences of numbers, plus improvements on the best-known recurrences known for other sequences.

\section{History and Outline}

The main source of the historical results is a paper \cite{FaasePaper} and website \cite{FaaseSite} by Faase, where the main motivation is to count the number of hamiltonian cycles in certain classes of graphs. Later on, in 2000, Desjarlais and Molina \cite{Desjarlais} discuss the number of spanning trees in $2 \times n$ and $3 \times n$ grid graphs. In 2004, Golin and Leung \cite{Golin} discuss a technique called \emph{unhooking} which will be used in this paper to reduce the problem of counting spanning trees in cylinder graphs to the problem of counting spanning trees in grid graphs.

In the first two papers, and this one, the general idea is the same: our goal is to count the number of spanning trees, but the method we use requires us to count other related objects, also. The paper by Faase appeals to the Transfer-Matrix Method, used widely in statistical mechanics (for more about the Transfer-Matrix Method, see [\cite{Stanley}]). The main distinction of this paper from \cite{Desjarlais} is the direct application of the Cayley-Hamilton Theorem to achieve recurrences for the sequences we are investigating. Overall, the results from this paper yield sequences for the number of spanning trees of the graphs $G \times P_n$ and $G \times C_n$ for any graph $G$. Along with these sequences, our methods find the minimal recurrence, generating function, and closed-form formulae for all of these sequences. As a consequence, we also find the sequences and recurrences for many, many other types of subgraphs.

The bulk of the paper focuses on the steps involved in finding the transition matrix for a given graph. In doing so, we will have to count other, related spanning forests with special properties.

\section{Notation.}

All of the graphs we will be dealing with depend on two parameters, which we will call $k$ and $n$. In all cases, we will think of $k$ as fixed and $n \to \infty$.

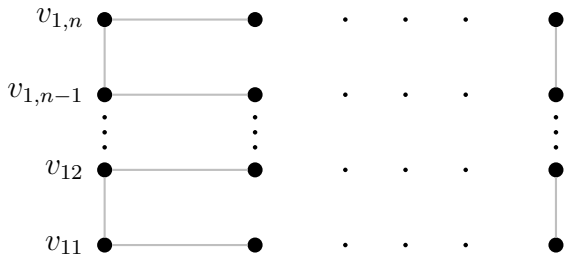
\begin{figure}[ht]
\begin{center}
\begin{pspicture}(0,0)(6,3)
\cnode*(0,0){0.1}{A}
\cnode*(0,1){0.1}{B}
\cnode*(0,2){0.1}{C}
\cnode*(0,3){0.1}{D}
\cnode*(2,0){0.1}{E}
\cnode*(2,1){0.1}{F}
\cnode*(2,2){0.1}{G}
\cnode*(2,3){0.1}{H}
\cnode*(6,0){0.1}{I}
\cnode*(6,1){0.1}{J}
\cnode*(6,2){0.1}{K}
\cnode*(6,3){0.1}{L}
\cnode*(0,1.3){0.03}{a}
\cnode*(0,1.5){0.03}{b}
\cnode*(0,1.7){0.03}{c}
\cnode*(2,1.3){0.03}{d}
\cnode*(2,1.5){0.03}{e}
\cnode*(2,1.7){0.03}{f}
\cnode*(6,1.3){0.03}{g}
\cnode*(6,1.5){0.03}{h}
\cnode*(6,1.7){0.03}{i}
\cnode*(3.2,0){0.03}{j}
\cnode*(4,0){0.03}{k}
\cnode*(4.8,0){0.03}{l}
\cnode*(3.2,1){0.03}{m}
\cnode*(4,1){0.03}{n}
\cnode*(4.8,1){0.03}{o}
\cnode*(3.2,2){0.03}{p}
\cnode*(4,2){0.03}{q}
\cnode*(4.8,2){0.03}{r}
\cnode*(3.2,3){0.03}{s}
\cnode*(4,3){0.03}{t}
\cnode*(4.8,3){0.03}{u}
\nput*{180}{A}{$v_{11}$}
\nput*{180}{B}{$v_{12}$}
\nput*{180}{C}{$v_{1,n-1}$}
\nput*{180}{D}{$v_{1,n}$}
\ncline[linecolor=lightgray]{-}{A}{B}
\ncline[linecolor=lightgray]{-}{C}{D}
\ncline[linecolor=lightgray]{-}{A}{E}
\ncline[linecolor=lightgray]{-}{B}{F}
\ncline[linecolor=lightgray]{-}{C}{G}
\ncline[linecolor=lightgray]{-}{D}{H}
\ncline[linecolor=lightgray]{-}{I}{J}
\ncline[linecolor=lightgray]{-}{K}{L}
\end{pspicture}
\caption{Labeling convention for grid graphs. \label{labeling}}
\end{center}
\end{figure}

\begin{definition}
The $k \times n$ \emph{grid graph} $G_k(n)$ is the simple graph with vertex and edge sets as follows:
\end{definition}
\begin{align*}
V\left(G_k(n)\right) &= \set{v_{ij} \mid 1 \leq i \leq k, 1 \leq j \leq n} \\
E\left(G_k(n)\right) &= \set{v_{i,j}v_{i',j'} \mid |i - i'| + |j - j'| = 1}
\end{align*}
In order to keep the diagrams clean, Figure \ref{labeling} shows the vertex naming conventions we will use.

\begin{figure}[ht]
\begin{center}
\begin{pspicture}(0,0)(2,2)
\cnode*(0,0){0.1}{A}
\cnode*(0,1){0.1}{B}
\cnode*(0,2){0.1}{C}
\cnode*(1,0){0.1}{D}
\cnode*(1,1){0.1}{E}
\cnode*(1,2){0.1}{F}
\cnode*(2,0){0.1}{G}
\cnode*(2,1){0.1}{H}
\cnode*(2,2){0.1}{I}
\ncline{-}{A}{B}
\ncline[linecolor=lightgray]{-}{B}{C}
\ncline{-}{A}{D}
\ncline[linecolor=lightgray]{-}{B}{E}
\ncline{-}{C}{F}
\ncline[linecolor=lightgray]{-}{D}{E}
\ncline{-}{E}{F}
\ncline[linecolor=lightgray]{-}{D}{G}
\ncline[linecolor=lightgray]{-}{E}{H}
\ncline{-}{F}{I}
\ncline{-}{G}{H}
\ncline{-}{H}{I}
\end{pspicture}
\caption{A forest in a $3 \times 3$ grid. \label{edge-ex}}
\end{center}
\end{figure}
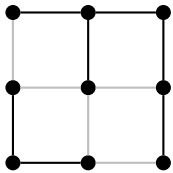

When showing examples, usually of spanning trees or spanning forests, we will always show the underlying graph in one form or another. A concrete example is given in figure \ref{edge-ex}: we will use black edges for edges in the subgraph exemplified; all unused edges will show up in light grey.

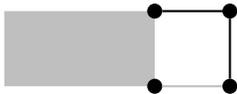
\begin{figure}[ht]
\begin{center}
\begin{pspicture}(0,0)(3,1)
\pspolygon*[linecolor=lightgray](0,0)(0,1)(2,1)(2,0)
\cnode*(2,1){.1}{A}
\cnode*(3,1){.1}{B}
\cnode*(3,0){.1}{C}
\cnode*(2,0){.1}{D}
\ncline{-}{A}{B}
\ncline{-}{B}{C}
\ncline[linecolor=lightgray]{-}{C}{D}
\end{pspicture}
\caption{An example of an arbitrary-sized graph with a specific end. \label{arb-ex}}
\end{center}
\end{figure}

When dealing with grids of arbitrary size, we will mainly be interested in the very right-most end of the grid, so we will represent the rest of the graph we do not care about by a gray box, as shown in figure \ref{arb-ex}.

The $k \times n$ \emph{cylinder graph} $C_k(n)$ can be obtained by ``wrapping'' the grid graph around, specifically by adding the following edges:
\[ E\left(C_k(n)\right) = E\left(G_k(n)\right) \bigcup \set{\set{v_{1,i},v_{n,i}} \mid 1 \leq i \leq k}. \]
Note that $C_k(n) = P_k \times C_n$.

The $k \times n$ \emph{torus graph} $T_k(n)$ can be obtained by ``wrapping'' the cylinder graph around the other way, specifically by adding the following edges:
\[ E\left(T_k(n)\right) = E\left(C_k(n)\right) \bigcup \set{\set{v_{i,1},v_{i,k}} \mid 1 \leq i \leq n} \]
Note that $T_k(n) = C_k \times C_n$.

%%%%%%%%%%%%%%%%
% We are excluding this stuff from the paper for now.
%%%%%%%%%%%%%%%%

%Finally, the \emph{generalized Star of David graph} $D_k(n)$ is the graph obtained by starting with a cycle $C_n$ and pasting a complete graph $K_k$ %along every edge of the cycle. For example, the Star of David graph is $D_3(6)$:

Throughout this paper, we will be dealing with partitions of the set $[k] = \set{1, 2, \ldots, k}$. We denote by $\mathcal{B}_k$ the set of all such partitions, and $B_k = |\mathcal{B}_k|$ are the Bell numbers. We will impose an ordering on $\mathcal{B}_k$, which we will call the lexicographic ordering on $\mathcal{B}_k$:

\begin{definition}
Given two partitions $P_1$ and $P_2$ of $[k]$, for $i \in [k]$, let $X_i$ be the block of $P_1$ containing $i$, and likewise $Y_i$
the block of $P_2$ containing $i$.  Let $j$ be the minimum value of $i$ such that $X_i \neq Y_i$.  Then $P_1 < P_2$ iff
\begin{enumerate}
\item $|P_1| < |P_2|$ or
\item $|P_1| = |P_2|$ and $X_j \prec Y_j$, where $\prec$ denotes normal lexicographic ordering.
\end{enumerate}
\end{definition}
For example, $\mathcal{B}_3$ in order is
\[ \mathcal{B}_3 = \set{\set{\set{1,2,3}}, \set{\set{1},\set{2,3}}, \set{\set{1,2},\set{3}}, \set{\set{1,3},\set{2}}, \set{\set{1},\set{2},\set{3}}} \]
However, we will use shorthand notation for set partitions as follows:
\[ \mathcal{B}_3 = \set{123,1/23,12/3,13/2,1/2/3}. \]
Since our examples will only deal with $k<10$, we will not have to worry about double-digit numbers on our shorthand notation.

We will find many recurrences in this paper, all pertaining to the number of spanning trees of the graphs mentioned above. Since we will be dealing with each type of graph separately, we will always denote by $T_n$ the number of spanning trees of whatever graph we are dealing with at the moment, which will be unambiguous.

\section{Grid Graphs: The Example For $k = 2$.}

What follows is mainly from \cite{Desjarlais} and is the inspiration for the other results on grid graphs. We would like to find a recurrence for $T_n$, which for now will represent the number of spanning trees in $G_2(n)$. If we started out with a spanning tree on $G_2(n-1)$, then there are three different ways to add the additional two vertices to still make a spanning tree on $G_2(n)$:
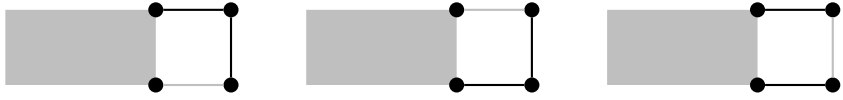
\begin{figure}[ht]
\begin{center}
\begin{pspicture}(0,0)(11,1)
\pspolygon*[linecolor=lightgray](0,0)(0,1)(2,1)(2,0)
\cnode*(2,1){.1}{A}
\cnode*(3,1){.1}{B}
\cnode*(3,0){.1}{C}
\cnode*(2,0){.1}{D}
\ncline{-}{A}{B}
\ncline{-}{B}{C}
\ncline[linecolor=lightgray]{-}{C}{D}
\pspolygon*[linecolor=lightgray](4,0)(4,1)(6,1)(6,0)
\cnode*(6,1){.1}{A}
\cnode*(7,1){.1}{B}
\cnode*(7,0){.1}{C}
\cnode*(6,0){.1}{D}
\ncline[linecolor=lightgray]{-}{A}{B}
\ncline{-}{C}{D}
\ncline{-}{B}{C}
\pspolygon*[linecolor=lightgray](8,0)(8,1)(10,1)(10,0)
\cnode*(10,1){.1}{A}
\cnode*(11,1){.1}{B}
\cnode*(11,0){.1}{C}
\cnode*(10,0){.1}{D}
\ncline[linecolor=lightgray]{-}{B}{C}
\ncline{-}{A}{B}
\ncline{-}{C}{D}
\end{pspicture}
\caption{Possible ways to extend a tree on $G_2(n-1)$ to obtain a tree on $G_2(n)$.}
\end{center}
\end{figure}

However, there is also a way to create a spanning tree on the $2 \times n$ grid from something that isn't a spanning tree on $G_k(n-1)$. Let $x = v_{1,n-1}$ and $y = v_{2,n-2}$ be the end vertices on $G_k(n-1)$. If we have a spanning forest on $G_k(n-1)$ with the property that there are two trees in the forest and $x$ and $y$ are in distinct trees, then we can append the following edges to create a spanning tree in $G_k(n)$:
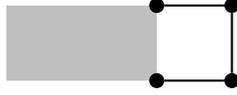
\begin{figure}[ht]
\begin{center}
\begin{pspicture}(0,0)(3,1)
\pspolygon*[linecolor=lightgray](0,0)(0,1)(2,1)(2,0)
\cnode*(2,1){.1}{A}
\cnode*(3,1){.1}{B}
\cnode*(3,0){.1}{C}
\cnode*(2,0){.1}{D}
\ncline{-}{A}{B}
\ncline{-}{B}{C}
\ncline{-}{C}{D}
\end{pspicture}
\caption{The only way to extend a certain forest on $G_2(n-1)$ to a tree on $G_2(n)$.}
\end{center}
\end{figure}

Therefore, in counting $T_n$ it is useful to also count $F_n$, which we define as the number of spanning forests in $G_k(n)$ consisting of two trees with the additional property that the end vertices $v_{1,n}$ and $v_{2,n}$ are in distinct trees. From the preceding two paragraphs we can now obtain the recurrence
\[ T_n = 3T_{n-1} + F_{n-1} \]
and through similar reasoning we can also find the recurrence
\[ F_n = 2T_{n-1} + F_{n-1} \]
At this point, let us note that we have enough information to find $T_n$ (or $F_n$) in time linear in $n$. However, our goal is to provide explicit recurrences for $T_n$ alone. If we let $v_n$ denote the column vector
\[ v_n = \left[ \begin{array}{c}
T_n \\
F_n
\end{array} \right] \]
And if we define the matrix $A$ by
\[ A = \left[\begin{array}{cc}
3 & 1 \\
2 & 1
\end{array}\right] \]
Then we satisfy
\[ Av_{n-1} = v_n. \]
With the starting conditions
\[ v_1 = \left[
\begin{array}{c}
1 \\
1
\end{array} \right]. \]
The characteristic polynomial of $A$ is
\[ \chi_\lambda(A) = \lambda^2 - 4\lambda + 1 \]
so by the Cayley-Hamilton Theorem, we satisfy
\[ A^2 - 4A + 1 = 0. \]
This can be re-written as
\[ A^2 = 4A - 1 \]
and if we multiply by the vector $v_n$ on the right we obtain
\[ \left[ \begin{array}{c} T_{n+2} \\ F_{n+2} \end{array} \right] = 4\left[ \begin{array}{c} T_{n+1} \\ F_{n+1} \end{array} \right] - \left[ \begin{array}{c} T_{n} \\ F_{n} \end{array} \right]. \]
Hence, we now see that $T_n$ and $F_n$ satisfy the same recurrence:
\begin{align*}
T_{n+2} &= 4T_{n+1} - T_n \\
F_{n+2} &= 4F_{n+1} - F_n
\end{align*}
with starting conditions
\[ \begin{array}{cc}
T_0 = 1 & T_1 = 4 \\
F_0 = 1 & F_1 = 3
\end{array}. \]
We now have all the information we need to obtain more information, such as the generating function and, finally, a closed-form formula for $T_n$. All of these items can be found in \cite{Desjarlais}.

\section{The General Case For Grid Graphs.}

We want to use the same ideas for general $k$, but it requires a bit more bookkeeping. To extend the idea of $F_n$ in the previous section, we need to consider partitions of $[k] = \set{1,2, \ldots, k}$ and the forests that come from these partitions.
\begin{definition}
Given a spanning forest $\mathcal{F}$ of $G_k(n)$, the \emph{partition induced by $\mathcal{F}$} is obtained from the equivalence relation
\[ i \sim j \iff v_{n,i}, v_{n,j} \mbox{ are in the same tree of } \mathcal{F}. \]
\end{definition}
For example, the partition induced by a spanning tree of $G_k(n)$ is $123\cdots n$ and the partition induced by the forest with no edges is $1/2/3/\cdots/n-1/n$.
\begin{definition}
Given a spanning forest $\mathcal{F}$ of $G_k(n)$ and a partition $P$ of $[k]$, we say that $\mathcal{F}$ \emph{is consistent with} $P$ if:
\begin{enumerate}
\item The number of trees in $\mathcal{F}$ is precisely $|P|$.
\item $P$ is the partition induced by $\mathcal{F}$.
\end{enumerate}
\end{definition}
\begin{definition}
Given a graph $G$ on $k$ vertices and a partition $P$ of $[k]$, let $T_G(P,n)$ be the number of spanning trees of the graph $G \times P_n$. We will often omit $G$ when it is clear from the context, or irrelevant. Recall that we have an ordering of partitions, so we will define $T_G(i,n) = T_G(P_i,n)$.
\end{definition}

In the previous section, since $B_2 = 2$, we were counting two things: $T_n$ , which corresponds to $T(12,n)$, and $F_n$, which corresponds to $T(1/2,n)$. Therefore, for arbitrary $k$ we are now tasked with counting $B_k$ different objects at once, so we are to find the $B_k \times B_k$ matrix that represents the $B_k$ simultaneous recurrences between these objects.

\begin{definition}
Define by $E_{n}$ the set of edges
\[ E_{n} = E(G_k(n)) \setminus E(G_k(n-1)) \]
Note that $\left|E_{n}\right| = 2k-1$ edges.
\end{definition}

Given some forest $\mathcal{F}$ of $G_k(n-1)$ and some subset $X \subseteq E_{n}$, we can combine the two to make a forest of $G_k(n)$. If we are only interested in the number of trees in the new forest and its induced partition, then we only need to know the same information from $\mathcal{F}$, and this is all independent of $n$. Therefore, we have the following definition:

\begin{definition}
Given two partitions $P_1$ and $P_2$ in $\mathcal{B}_k$, a subset $X \subseteq E_{n}$ \emph{transfers from $P_1$ to $P_2$} if a forest consistent with $P_1$ becomes a forest consistent with $P_2$ after the addition of $X$.
\end{definition}
\begin{example}
Figure \ref{4x4-ex} shows a spanning forest of $G_4(4)$ where, from left to right, the edges transfer from $1/23/4$ to $1234$, from $1234$ to $12/34$, and from $12/34$ to $1/2/34$.
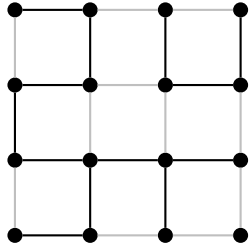
\begin{figure}[ht]
\begin{center}
\begin{pspicture}(0,0)(2,3)
\cnode*(0,0){0.1}{A}
\cnode*(0,1){0.1}{B}
\cnode*(0,2){0.1}{C}
\cnode*(0,3){0.1}{D}
\cnode*(1,0){0.1}{E}
\cnode*(1,1){0.1}{F}
\cnode*(1,2){0.1}{G}
\cnode*(1,3){0.1}{H}
\cnode*(2,0){0.1}{I}
\cnode*(2,1){0.1}{J}
\cnode*(2,2){0.1}{K}
\cnode*(2,3){0.1}{L}
\cnode*(3,0){0.1}{M}
\cnode*(3,1){0.1}{N}
\cnode*(3,2){0.1}{O}
\cnode*(3,3){0.1}{P}
\ncline[linecolor=lightgray]{-}{A}{B}
\ncline{-}{B}{C}
\ncline[linecolor=lightgray]{-}{C}{D}
\ncline{-}{A}{E}
\ncline{-}{B}{F}
\ncline{-}{C}{G}
\ncline{-}{D}{H}
\ncline{-}{E}{F}
\ncline[linecolor=lightgray]{-}{F}{G}
\ncline{-}{G}{H}
\ncline[linecolor=lightgray]{-}{E}{I}
\ncline{-}{F}{J}
\ncline[linecolor=lightgray]{-}{G}{K}
\ncline[linecolor=lightgray]{-}{H}{L}
\ncline{-}{I}{J}
\ncline[linecolor=lightgray]{-}{J}{K}
\ncline{-}{K}{L}
\ncline[linecolor=lightgray]{-}{I}{M}
\ncline{-}{J}{N}
\ncline{-}{K}{O}
\ncline[linecolor=lightgray]{-}{L}{P}
\ncline[linecolor=lightgray]{-}{M}{N}
\ncline[linecolor=lightgray]{-}{N}{O}
\ncline{-}{O}{P}
\end{pspicture}
\caption{An example of a spanning forest of $G_4(4)$. \label{4x4-ex}}
\end{center}
\end{figure}
\end{example}

Therefore, we can define the $B_k \times B_k$ matrix $A_k$ by the following:
\[ A_k(i,j) = |\set{A \subseteq E_{n+1} \mid A \mbox{ is compatible from } P_j \mbox{ to } P_i}|. \]
The $2 \times 2$ matrix in the previous section is $A_2$. Brute-force search with straightforward Mathematica code \cite{RaffSite} can produce more matrices:

\[ A_3 = \begin{bmatrix}
8 & 3 & 3 & 4 & 1 \\
4 & 3 & 2 & 2 & 1 \\
4 & 2 & 3 & 2 & 1 \\
1 & 0 & 0 & 1 & 0 \\
3 & 2 & 2 & 2 & 1
\end{bmatrix}
\]

\[ A_4 =\left[ \begin{array}{p{6pt}p{6pt}p{6pt}p{6pt}p{6pt}p{6pt}p{6pt}p{6pt}p{6pt}p{6pt}p{6pt}p{6pt}p{6pt}p{6pt}p{6pt}}
21 & 8 & 9 & 11 & 8 & 14 & 11 & 15 & 3 & 3 & 4 & 3 & 4 & 5 & 1 \\
9 & 8 & 6 & 4 & 4 & 6 & 5 & 8 & 3 & 3 & 4 & 2 & 2 & 2 & 1 \\
6 & 4 & 9 & 4 & 4 & 4 & 4 & 4 & 3 & 2 & 2 & 3 & 2 & 2 & 1 \\
3 & 0 & 0 & 3 & 1 & 2 & 1 & 2 & 0 & 0 & 0 & 0 & 1 & 1 & 0 \\
9 & 4 & 6 & 5 & 8 & 6 & 4 & 8 & 2 & 3 & 2 & 3 & 4 & 2 & 1 \\
1 & 0 & 0 & 1 & 0 & 3 & 1 & 0 & 0 & 0 & 0 & 0 & 0 & 1 & 0 \\
3 & 1 & 0 & 1 & 0 & 2 & 3 & 2 & 0 & 0 & 1 & 0 & 0 & 1 & 0 \\
0 & 0 & 0 & 0 & 0 & 0 & 0 & 1 & 0 & 0 & 0 & 0 & 0 & 0 & 0 \\
5 & 4 & 6 & 4 & 3 & 4 & 3 & 4 & 3 & 2 & 2 & 2 & 2 & 2 & 1 \\
5 & 4 & 4 & 3 & 4 & 6 & 3 & 4 & 2 & 3 & 2 & 2 & 2 & 2 & 1 \\
1 & 1 & 0 & 0 & 0 & 0 & 1 & 2 & 0 & 0 & 1 & 0 & 0 & 0 & 0 \\
5 & 3 & 6 & 3 & 4 & 4 & 4 & 4 & 2 & 2 & 2 & 3 & 2 & 2 & 1 \\
1 & 0 & 0 & 1 & 1 & 0 & 0 & 2 & 0 & 0 & 0 & 0 & 1 & 0 & 0 \\
1 & 0 & 0 & 1 & 0 & 2 & 1 & 0 & 0 & 0 & 0 & 0 & 0 & 1 & 0 \\
4 & 3 & 4 & 3 & 3 & 4 & 3 & 4 & 2 & 2 & 2 & 2 & 2 & 2 & 1
\end{array} \right]
\]

$A_5,A_6,$ and $A_7$ have also been found; they are shown in \cite{RaffSite}. Once these matrices are known, then everything about the sequence of spanning trees can be found. The following table shows some results obtained for grid graphs; results obtained for arbitrary graphs of the form $G \times P_n$ for all graphs $G$ with at most five vertices are in \cite{RaffSite}.

\renewcommand{\arraystretch}{1.1}

\begin{tabular}{|l|}
\hline
$G_2(n):$ (\cite{Desjarlais}) \\
$T_n = 4T_{n-1}-T_{n-2}$ \\
Sequence: $\set{1,4,15,56,209, \ldots}$ (OEIS A001353)\\
Generating Function: $\frac{x}{1 - 4x + x^2}$\\
\hline
$G_3(n):$ (\cite{FaaseSite}) \\
$T_n = 15T_{n-1} - 32T_{n-2} + 15T_{n-3} - T_{n-4}$ \\
Sequence: $\set{1, 15, 192, 2415, 30305, \ldots}$ (OEIS A006238)\\
Generating Function: $\frac{3x(1 + 49x + 1152x^2)}{1 + 24x - 24x^2 + x^3}$ \\
\hline
$G_4(n):$ (\cite{FaaseSite}) \\
$T_n = 56T_{n-1} - 672T_{n-2} + 2632T_{n-3} - 4094T_{n-4} + 2632T_{n-5} - 672T_{n-6} + 56T_{n-7} - T_{n-8}$ \\
Sequence: $\set{1, 56, 2415, 100352, 4140081, \ldots}$ (OEIS A003696)\\
Generating Function: $\frac{16x(1 + 12x + x^2)}{1 - 204x + 1190x^2 - 204x^3 + x^4}$ \\
\hline
$G_5(n):$ (\cite{FaaseSite}, with improvements from this paper) \\
$T_n = 209T_{n-1} - 11936T_{n-2} + 274208T_{n-3} - 3112032T_{n-4} + 19456019T_{n-5}$ \\
$- 70651107T_{n-6} + 152325888T_{n-7} - 196664896T_{n-8} + 152325888T_{n-9}$ \\
$- 70651107T_{n-10} + 19456019T_{n-11} - 3112032T_{n-12} + 274208T_{n-13}$ \\
$- 11936T_{n-14} + 209T_{n-15} - T_{n-16}$ \\
Sequence: $\set{1, 209, 30305, 4140081, 557568000, \ldots}$ (OEIS A003779)  \\
Generating Function: $\frac{125x(1 + 4656x + 10616686x^2 + 23432228161x^3 + 51714958501250x^4)}{1 + 2255x - 105985x^2 + 105985x^3 - 2255x^4 + x^5}$\\
\hline
$G_6(n):$ (new) \\
$T_n = 780 T_{n-1}-194881 T_{n-2}+22377420 T_{n-3}-1419219792 T_{n-4}$ \\
           $+ 55284715980T_{n-5}-1410775106597 T_{n-6}+24574215822780 T_{n-7} $\\
           $-300429297446885 T_{n-8}+2629946465331120 T_{n-9}-16741727755133760 T_{n-10} $ \\
           $+78475174345180080 T_{n-11}-273689714665707178 T_{n-12}+716370537293731320 T_{n-13} $ \\
           $-1417056251105102122 T_{n-14}+2129255507292156360 T_{n-15}-2437932520099475424 T_{n-16} $ \\
           $+2129255507292156360 T_{n-17}-1417056251105102122 T_{n-18}+716370537293731320 T_{n-19} $ \\
           $-273689714665707178 T_{n-20}+78475174345180080 T_{n-21}-16741727755133760 T_{n-22} $ \\
           $+2629946465331120 T_{n-23}-300429297446885 T_{n-24}+24574215822780 T_{n-25} $ \\
           $-1410775106597 T_{n-26} +55284715980 T_{n-27}-1419219792 T_{n-28}+22377420 T_{n-29} $ \\
           $-194881 T_{n-30}+780 T_{n-31}-T_{n-32} $ \\
Sequence: $\set{1, 780, 380160, 170537640, 74795194705, \ldots}$ (OEIS A139400) \\
Generating Function: See \cite{RaffSite} \\
\hline
\end{tabular}

\section{Extending to Generalized Graphs of the Form $G \times P_n$}

For the results above, it was not necessary that the graph we were dealing with was a grid. We could have repeated the same process as above for any sequences of graphs $G_n$ defined by
\[ G_n = G \times P_n \]
for some predefined graph $G$. In fact, the Mathematica code in the appendix handles any such general case. Therefore, it leads to the following theorem:

\begin{theorem}
Let a graph $G$ be given with $k$ vertices, and define the sequence of graphs $\set{G_n}$ by $G_n = G \times P_n$. Then there is a $B_k \times B_k$ matrix $M$ and a vector $v$, both taking on integer values, such that
\[ T_n = M^nv[1] \]
where $T_n$ is the number of spanning trees in $G_n$. Furthermore, $M^nv[i]$ lists the number of spanning forests consistent with $P_i$ in $G_n$.
\end{theorem}

\begin{corollary}
Let a graph $G$ be given with $k$ vertices, and consider the sequence $\set{T_n}$. Then $T_n$ satisfies a linear recurrence of order $B_k$.
\end{corollary}

From investigations, we have a few conjectures:

\begin{conjecture}
For the matrix $M$ given in the theorem above, the characteristic polynomial $\chi_\lambda(M)$ factors over the integers into monomials whose degree is always a power of 2.
\end{conjecture}

\begin{conjecture}
For any graph $G$, the recurrence $\set{T_n}$ satisfies a linear recurrence whose coefficients alternate in sign.
\end{conjecture}

\begin{conjecture}
The recurrence for the grid graph $G_k(n)$ has order $2^{k-1}$.
\end{conjecture}

\begin{conjecture}
The recurrence for the graph $K_k \times P_n$ has order $k$.
\end{conjecture}

For the time being, we will only prove the special case of Conjecture 3 for the grid graphs $G_2(n)$. We will give a combinatorial proof that we hope can be adjusted accordingly to the higher cases. To aid in the proof, we will introduce the concept of \emph{grid addition}, which is simply a shorthand way of creating the union of two grids.

\begin{definition}
If $G_1$ is a $k \times n_1$ grid and $G_2$ is a $k \times n_2$ grid, then $G_1 + G_2$ is the $k \times (n_1 + n_2 - 1)$ grid defined as the graph obtained by identifying the right-most vertices of $G_1$ with the left-most vertices of $G_2$. Any overlapping edges remain.
\end{definition}

\begin{example}
Figure \ref{grid-add} shows the addition of a tree on $G_2(3)$ with a tree on $G_2(2)$ to obtain a subgraph of $G_2(4)$.
\begin{figure}[ht]
\begin{center}
\begin{pspicture}(0,0)(10,1)
\cnode*(0,0){.1}{A}
\cnode*(0,1){.1}{B}
\cnode*(1,0){.1}{C}
\cnode*(1,1){.1}{D}
\cnode*(2,0){.1}{E}
\cnode*(2,1){.1}{F}
\rput{0}(3,0.5){$+$}
\cnode*(4,0){.1}{G}
\cnode*(4,1){.1}{H}
\cnode*(5,0){.1}{I}
\cnode*(5,1){.1}{J}
\ncline{-}{A}{C}
\ncline{-}{B}{D}
\ncline{-}{C}{D}
\ncline{-}{D}{F}
\ncline{-}{C}{E}
\ncline{-}{G}{I}
\ncline{-}{H}{J}
\ncline{-}{I}{J}
\rput{0}(6,0.5){$=$}
\cnode*(7,0){.1}{A}
\cnode*(7,1){.1}{B}
\cnode*(8,0){.1}{C}
\cnode*(8,1){.1}{D}
\cnode*(9,0){.1}{E}
\cnode*(9,1){.1}{F}
\cnode*(9,0){.1}{G}
\cnode*(9,1){.1}{H}
\cnode*(10,0){.1}{I}
\cnode*(10,1){.1}{J}
\ncline{-}{A}{C}
\ncline{-}{B}{D}
\ncline{-}{C}{D}
\ncline{-}{D}{F}
\ncline{-}{C}{E}
\ncline{-}{G}{I}
\ncline{-}{H}{J}
\ncline{-}{I}{J}
\end{pspicture}
\caption{An example of grid addition. \label{grid-add}}
\end{center}
\end{figure}
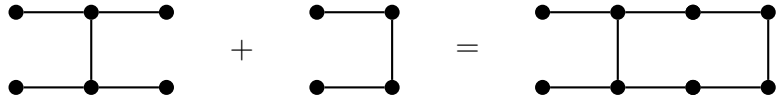
\end{example}

\begin{theorem}
The number of spanning trees of the graphs $G_2(n)$ satisfies the linear recurrence $T_n = 4T_{n-1} - T_{n-2}$ with the initial conditions $T_1 = 1$, $T_2 = 4$.
\end{theorem}
\begin{proof}
Showing the initial conditions is a minor exercise. We will prove this recurrence in the equivalent form $T_n + T_{n-2} = 4T_{n-1}$. Let $\mathcal{T}_k$ denote the set of spanning trees of the graph $G_2(k)$. We will associate $T_{n-2}$ with the set $\mathcal{T}_{n-2}$ with an addition at the end, as shown by Figure \ref{int-T-n-2}.
\begin{figure}[ht]
\begin{center}
\begin{pspicture}(6,1)
\pspolygon*[linecolor=lightgray](0,0)(0,1)(4,1)(4,0)
\cnode*(4,0){.1}{A}
\cnode*(4,1){.1}{B}
\cnode*(5,0){.1}{C}
\cnode*(5,1){.1}{D}
\cnode*(6,0){.1}{E}
\cnode*(6,1){.1}{F}
\ncline{-}{A}{C}
\ncline{-}{C}{E}
\ncline{-}{B}{D}
\ncline{D}{F}
\end{pspicture}
\caption{How we interpret $T_{n-2}$. \label{int-T-n-2}}
\end{center}
\end{figure}
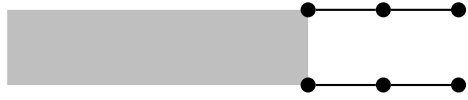
In this way, we can think of $\mathcal{T}_{n-2}$ as being trees of $G_2(n)$. Similarly, as Figure \ref{int-4T-n-1} shows, we will associate $4T_{n-1}$ with the set of trees from $\mathcal{T}_{n-1}$ with each of the four trees of $G_2(2)$ added at the end.
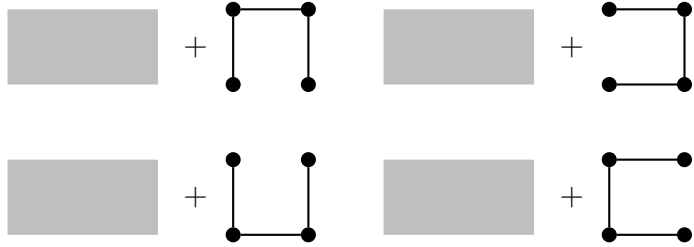
\begin{figure}[ht]
\begin{center}
\begin{pspicture}(9,3)
\pspolygon*[linecolor=lightgray](0,0)(0,1)(2,1)(2,0)
\rput{0}(2.5,0.5){$+$}
\cnode*(3,0){.1}{A}
\cnode*(3,1){.1}{B}
\cnode*(4,0){.1}{C}
\cnode*(4,1){.1}{D}
\ncline{-}{A}{B}
\ncline{-}{A}{C}
\ncline{-}{C}{D}
\pspolygon*[linecolor=lightgray](5,0)(5,1)(7,1)(7,0)
\rput{0}(7.5,0.5){$+$}
\cnode*(8,0){.1}{A}
\cnode*(8,1){.1}{B}
\cnode*(9,0){.1}{C}
\cnode*(9,1){.1}{D}
\ncline{-}{A}{B}
\ncline{-}{B}{D}
\ncline{-}{A}{C}
\pspolygon*[linecolor=lightgray](0,2)(0,3)(2,3)(2,2)
\rput{0}(2.5,2.5){$+$}
\cnode*(3,2){.1}{A}
\cnode*(3,3){.1}{B}
\cnode*(4,2){.1}{C}
\cnode*(4,3){.1}{D}
\ncline{-}{B}{D}
\ncline{-}{A}{B}
\ncline{-}{C}{D}
\pspolygon*[linecolor=lightgray](5,2)(5,3)(7,3)(7,2)
\rput{0}(7.5,2.5){$+$}
\cnode*(8,2){.1}{A}
\cnode*(8,3){.1}{B}
\cnode*(9,2){.1}{C}
\cnode*(9,3){.1}{D}
\ncline{-}{A}{C}
\ncline{-}{B}{D}
\ncline{-}{C}{D}
\end{pspicture}
\caption{How we interpret $4T_{n-1}$. \label{int-4T-n-1}}
\end{center}
\end{figure}
If we have a tree from $\mathcal{T}_n$, then we can decompose it depending on what the ending of the tree looks like. Figure \ref{decomp-Tn} shows all of the possibilities, along with their decompositions. Note that the decompositions are of the same form as we dictated for $4T_{n-1}$.
\begin{figure}[ht]
\begin{center}
\begin{pspicture}(10,13)
\pspolygon*[linecolor=lightgray](0,0)(0,1)(2,1)(2,0)
\cnode*(2,0){.1}{A}
\cnode*(2,1){.1}{B}
\cnode*(3,0){.1}{C}
\cnode*(3,1){.1}{D}
\ncline{-}{A}{C}
\ncline{-}{B}{D}
\rput{0}(4.5,0.5){$\to$}
\pspolygon*[linecolor=lightgray](6,0)(6,1)(8,1)(8,0)
\rput{0}(8.5,0.5){$+$}
\cnode*(8,0){.1}{E1}
\cnode*(8,1){.1}{E2}
\cnode*(9,0){.1}{A}
\cnode*(9,1){.1}{B}
\cnode*(10,0){.1}{C}
\cnode*(10,1){.1}{D}
\ncline{-}{A}{C}
\ncline{-}{A}{B}
\ncline{-}{B}{D}
\pspolygon*[linecolor=lightgray](0,2)(0,3)(2,3)(2,2)
\cnode*(2,2){.1}{A}
\cnode*(2,3){.1}{B}
\cnode*(3,2){.1}{C}
\cnode*(3,3){.1}{D}
\ncline{-}{A}{C}
\ncline{-}{B}{D}
\ncline{-}{A}{B}
\rput{0}(4.5,2.5){$\to$}
\pspolygon*[linecolor=lightgray](6,2)(6,3)(8,3)(8,2)
\rput{0}(8.5,2.5){$+$}
\cnode*(8,2){.1}{E1}
\cnode*(8,3){.1}{E2}
\cnode*(9,2){.1}{A}
\cnode*(9,3){.1}{B}
\cnode*(10,2){.1}{C}
\cnode*(10,3){.1}{D}
\ncline{-}{E1}{E2}
\ncline{-}{A}{C}
\ncline{-}{A}{B}
\ncline{-}{B}{D}
\pspolygon*[linecolor=lightgray](0,4)(0,5)(2,5)(2,4)
\cnode*(2,4){.1}{A}
\cnode*(2,5){.1}{B}
\cnode*(3,4){.1}{C}
\cnode*(3,5){.1}{D}
\ncline{-}{B}{D}
\ncline{-}{C}{D}
\rput{0}(4.5,4.5){$\to$}
\pspolygon*[linecolor=lightgray](6,4)(6,5)(8,5)(8,4)
\rput{0}(8.5,4.5){$+$}
\cnode*(8,4){.1}{E1}
\cnode*(8,5){.1}{E2}
\cnode*(9,4){.1}{A}
\cnode*(9,5){.1}{B}
\cnode*(10,4){.1}{C}
\cnode*(10,5){.1}{D}
\ncline{-}{B}{D}
\ncline{-}{A}{B}
\ncline{-}{C}{D}
\pspolygon*[linecolor=lightgray](0,6)(0,7)(2,7)(2,6)
\cnode*(2,6){.1}{A}
\cnode*(2,7){.1}{B}
\cnode*(3,6){.1}{C}
\cnode*(3,7){.1}{D}
\ncline{-}{A}{B}
\ncline{-}{B}{D}
\ncline{-}{C}{D}
\rput{0}(4.5,6.5){$\to$}
\pspolygon*[linecolor=lightgray](6,6)(6,7)(8,7)(8,6)
\rput{0}(8.5,6.5){$+$}
\cnode*(8,6){.1}{E1}
\cnode*(8,7){.1}{E2}
\cnode*(9,6){.1}{A}
\cnode*(9,7){.1}{B}
\cnode*(10,6){.1}{C}
\cnode*(10,7){.1}{D}
\ncline{-}{E1}{E2}
\ncline{-}{A}{B}
\ncline{-}{B}{D}
\ncline{-}{C}{D}
\pspolygon*[linecolor=lightgray](0,8)(0,9)(2,9)(2,8)
\cnode*(2,8){.1}{A}
\cnode*(2,9){.1}{B}
\cnode*(3,8){.1}{C}
\cnode*(3,9){.1}{D}
\ncline{-}{A}{C}
\ncline{-}{C}{D}
\rput{0}(4.5,8.5){$\to$}
\pspolygon*[linecolor=lightgray](6,8)(6,9)(8,9)(8,8)
\rput{0}(8.5,8.5){$+$}
\cnode*(8,8){.1}{E1}
\cnode*(8,9){.1}{E2}
\cnode*(9,8){.1}{A}
\cnode*(9,9){.1}{B}
\cnode*(10,8){.1}{C}
\cnode*(10,9){.1}{D}
\ncline{-}{A}{B}
\ncline{-}{A}{C}
\ncline{-}{C}{D}
\pspolygon*[linecolor=lightgray](0,10)(0,11)(2,11)(2,10)
\cnode*(2,10){.1}{A}
\cnode*(2,11){.1}{B}
\cnode*(3,10){.1}{C}
\cnode*(3,11){.1}{D}
\ncline{-}{A}{B}
\ncline{-}{A}{C}
\ncline{-}{C}{D}
\rput{0}(4.5,10.5){$\to$}
\pspolygon*[linecolor=lightgray](6,10)(6,11)(8,11)(8,10)
\rput{0}(8.5,10.5){$+$}
\cnode*(8,10){.1}{E1}
\cnode*(8,11){.1}{E2}
\cnode*(9,10){.1}{A}
\cnode*(9,11){.1}{B}
\cnode*(10,10){.1}{C}
\cnode*(10,11){.1}{D}
\ncline{-}{E1}{E2}
\ncline{-}{A}{B}
\ncline{-}{A}{C}
\ncline{-}{C}{D}
\pspolygon*[linecolor=lightgray](0,12)(0,13)(2,13)(2,12)
\cnode*(2,12){.1}{A}
\cnode*(2,13){.1}{B}
\cnode*(3,12){.1}{C}
\cnode*(3,13){.1}{D}
\ncline{-}{A}{C}
\ncline{-}{B}{D}
\ncline{-}{C}{D}
\rput{0}(4.5,12.5){$\to$}
\pspolygon*[linecolor=lightgray](6,12)(6,13)(8,13)(8,12)
\rput{0}(8.5,12.5){$+$}
\cnode*(8,12){.1}{E1}
\cnode*(8,13){.1}{E2}
\cnode*(9,12){.1}{A}
\cnode*(9,13){.1}{B}
\cnode*(10,12){.1}{C}
\cnode*(10,13){.1}{D}
\ncline{-}{E1}{E2}
\ncline{-}{A}{C}
\ncline{-}{B}{D}
\ncline{-}{C}{D}
\end{pspicture}
\caption{Endings and decompositions for elements of $\mathcal{T}_n$. \label{decomp-Tn}}
\end{center}
\end{figure}
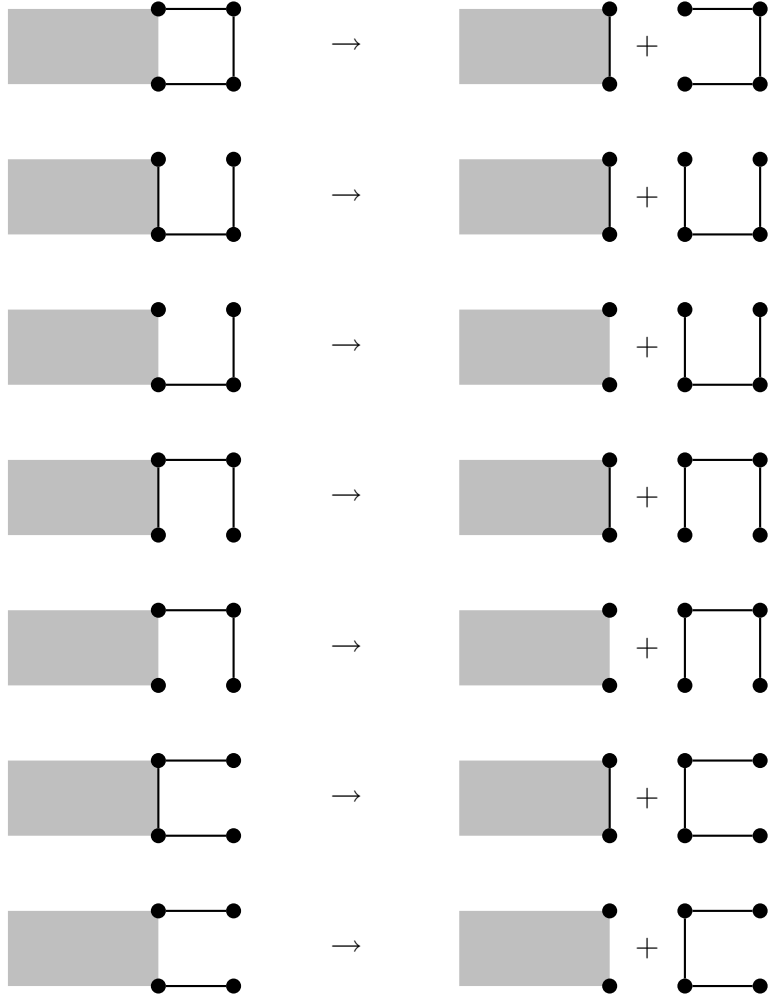
Similarly, if we have a tree from $\mathcal{T}_{n-2}$ modified as explained above, then Figure \ref{decomp-Tn-2} shows the decomposition. Again, note that the decompositions are of the same form as we dictated for $4T_{n-1}$.
\begin{figure}[ht]
\begin{center}
\begin{pspicture}(12,1)
\pspolygon*[linecolor=lightgray](0,0)(0,1)(2,1)(2,0)
\cnode*(2,0){.1}{A}
\cnode*(2,1){.1}{B}
\cnode*(3,0){.1}{C}
\cnode*(3,1){.1}{D}
\cnode*(4,0){.1}{E}
\cnode*(4,1){.1}{F}
\ncline{-}{A}{C}
\ncline{-}{B}{D}
\ncline{-}{C}{E}
\ncline{-}{D}{F}
\rput{0}(5.5,0.5){$\to$}
\pspolygon*[linecolor=lightgray](7,0)(7,1)(9,1)(9,0)
\rput{0}(10.5,0.5){$+$}
\cnode*(9,0){.1}{A}
\cnode*(9,1){.1}{B}
\cnode*(10,0){.1}{E1}
\cnode*(10,1){.1}{E2}
\cnode*(11,0){.1}{C}
\cnode*(11,1){.1}{D}
\cnode*(12,0){.1}{E}
\cnode*(12,1){.1}{F}
\ncline{-}{A}{E1}
\ncline{-}{B}{E2}
\ncline{-}{C}{E}
\ncline{-}{D}{F}
\ncline{-}{E}{F}
\end{pspicture}
\caption{Ending and decomposition for elements of $\mathcal{T}_{n-2}$. \label{decomp-Tn-2}}
\end{center}
\end{figure}
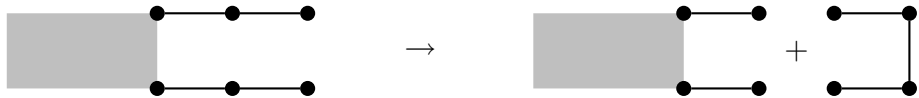
The reader can verify that the map described is invertible, yielding the desired bijection.
\end{proof}

\section{Extending to Cylinder Graphs}

\begin{figure}[ht]
\begin{center}
\begin{pspicture}(0,0)(2,2)
\psline[linearc=.1,linecolor=lightgray]{-}(2,2)(2.25,2.25)(-.25,2.25)(0,2)
\psline[linearc=.1,linecolor=lightgray]{-}(2,1)(2.25,1.25)(-.25,1.25)(0,1)
\psline[linearc=.1,linecolor=lightgray]{-}(2,0)(2.25,0.25)(-.25,0.25)(0,0)
\cnode*(0,0){0.1}{A}
\cnode*(0,1){0.1}{B}
\cnode*(0,2){0.1}{C}
\cnode*(1,0){0.1}{D}
\cnode*(1,1){0.1}{E}
\cnode*(1,2){0.1}{F}
\cnode*(2,0){0.1}{G}
\cnode*(2,1){0.1}{H}
\cnode*(2,2){0.1}{I}
\ncline{-}{A}{B}
\ncline[linecolor=lightgray]{-}{B}{C}
\ncline{-}{A}{D}
\ncline[linecolor=lightgray]{-}{B}{E}
\ncline{-}{C}{F}
\ncline[linecolor=lightgray]{-}{D}{E}
\ncline{-}{E}{F}
\ncline[linecolor=lightgray]{-}{D}{G}
\ncline[linecolor=lightgray]{-}{E}{H}
\ncline{-}{F}{I}
\ncline{-}{G}{H}
\ncline{-}{H}{I}
\end{pspicture}
\caption{Example for cylinder. \label{cyl-ex}}
\end{center}
\end{figure}
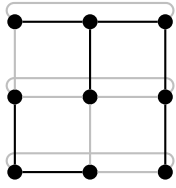

In this section we will discuss the changes necessary to extend the above arguments to find recurrences for cylinder graphs and generalized cylinder graphs. We shall take advantage of the ``unhooking'' technique covered in \cite{Golin}. The technique is a reduction from a cylinder graph to a grid graph. Recall that the vertex sets of $C_k(n)$ and $G_k(n)$ are the same.
\begin{definition}
For a given $k$, we define $\mathcal{E}_k$ by
\[ \mathcal{E}_k = E(C_k(n))\setminus E(G_k(n)) \]
\end{definition}
If we unhook (i.e. remove) the edges in $\mathcal{E}_k$ then what we have left is precisely $G_k(n)$. Now we have to consider what structures in $G_k(n)$ yield a spanning tree in $C_k(n)$ by the addition of some subset of edges from $\mathcal{E}_k$. Since we are going to add edges that go from one end of the grid to another, we must look at both ends of the grid now, as opposed to only looking at one end. For example, Figure \ref{cyl-ex} shows a spanning forest of $G_3(3)$ will never yield a spanning tree of $G_3(n)$ for any $n > 3$ through the method described in the previous sections, but this spanning forest would create two different spanning trees of $C_3(3)$ through the addition of either edge $v_{1,1}v_{3,1}$ or $v_{1,2}v_{3,2}$.

Therefore, we can keep the same basic idea used with grid graphs, with some modifications. We must now keep track of how our spanning forests affects the vertices \emph{at each end}.
\begin{definition}
Given a spanning forest $\mathcal{F}$ of $G_k(n)$, the partition $P$ of $[2k]$ \emph{induced by} $\mathcal{F}$ is obtained from the equivalence relation
\[ i \sim j \iff v_i, v_j \mbox{ are in the same tree of } \mathcal{F} \]
where we identify the vertices $v_1, v_2, \ldots , v_k$ with $v_{1,1}, v_{1,2}, \ldots , v_{1,k}$, respectively, and the vertices $v_{k+1},v_{k+2},\ldots,v_{2k}$ with $v_{n,1}, v_{n,2}, \ldots, v_{n,k}$, respectively.
\end{definition}
\begin{definition}
Given a spanning forest $\mathcal{F}$ of $G_k(n)$ and a partition $P$ of $[2k]$, we say that $\mathcal{F}$ \emph{is cylindrically consistent with} $P$ if:
\begin{enumerate}
\item The number of trees in $\mathcal{F}$ is precisely $|P|$.
\item $P$ is the partition induced by $\mathcal{F}$.
\end{enumerate}
\end{definition}
For example, the forest shown in Figure \ref{cyl-ex} is consistent with the partition $12/3456$. It's important to know what partition a certain forest of $G_k(n)$ is cylindrically consistent with, as that determines how many different ways edges can be added to achieve a spanning tree of $C_k(n)$. Since each spanning tree of $C_k(n)$ is uniquely determined by the underlying spanning forest of $G_k(n)$ and the extra edges from $\mathcal{E}_k$, we have all the information we need to count the number of spanning trees of $C_k(n)$.
\begin{definition}
For a given $k$, the \emph{tree-counting vector} $d_k$ is the vector, indexed by the partitions of $[2k]$, such that $d_k(i)$ is the number of ways that edges from $E(C_k(n)) \setminus E(G_k(n))$ can be added to get from a forest cylindrically consistent with partition $i$ to a spanning tree of $C_k(n)$. Notice that this is independent of $n$.
\end{definition}
It can be verified that the following information produces $d_2$:
\[\begin{array}{|c|c|}
\hline
1234 & 1 \\
\hline
1/234 & 1 \\
\hline
12/34 & 2 \\
\hline
134/2 & 1 \\
\hline
123/4 & 1 \\
\hline
14/23 & 2 \\
\hline
124/3 & 1 \\
\hline
13/24 & 0 \\
\hline
1/2/34 & 1 \\
\hline
1/23/4 & 1 \\
\hline
1/24/3 & 0 \\
\hline
12/3/4 & 1 \\
\hline
13/2/4 & 0 \\
\hline
14/2/3 & 1 \\
\hline
1/2/3/4 & 0 \\
\hline
\end{array}
\]
\[ d_2 = (1,1,2,1,1,2,1,0,1,1,0,1,0,1,0) \]
To count the number of spanning trees for $C_k(n)$ we can produce the $B_{2k} \times B_{2k}$ matrix in the same way as we did for the grid graphs, and using this matrix we can find the number of spanning forests of $G_k(n)$ consistent with each of the partitions of $\mathcal{B}_{2k}$, which can be expressed as a vector of length $B_{2k}$. Then, when we take the dot product of this vector with $d_k$, we obtain the number of spanning trees of $C_k(n)$. For example, it can be verified that the following is the matrix related to $C_2(n)$:
\begin{equation*}
A =
\left[\begin{array}{p{6pt}p{6pt}p{6pt}p{6pt}p{6pt}p{6pt}p{6pt}p{6pt}p{6pt}p{6pt}p{6pt}p{6pt}p{6pt}p{6pt}p{6pt}}
3 & 0 & 0 & 0 & 1 & 1 & 1 & 1 & 0 & 0 & 0 & 0 & 0 & 0 & 0 \\
 0 & 3 & 0 & 0 & 0 & 1 & 0 & 1 & 0 & 1 & 1 & 0 & 0 & 0 & 0 \\
 1 & 0 & 3 & 0 & 1 & 0 & 1 & 0 & 0 & 0 & 0 & 1 & 0 & 0 & 0 \\
 0 & 0 & 0 & 3 & 0 & 1 & 0 & 1 & 0 & 0 & 0 & 0 & 1 & 1 & 0 \\
 1 & 0 & 0 & 0 & 1 & 0 & 0 & 0 & 0 & 0 & 0 & 0 & 0 & 0 & 0 \\
 0 & 0 & 0 & 0 & 0 & 1 & 0 & 0 & 0 & 0 & 0 & 0 & 0 & 0 & 0 \\
 1 & 0 & 0 & 0 & 0 & 0 & 1 & 0 & 0 & 0 & 0 & 0 & 0 & 0 & 0 \\
 0 & 0 & 0 & 0 & 0 & 0 & 0 & 1 & 0 & 0 & 0 & 0 & 0 & 0 & 0 \\
 0 & 1 & 0 & 1 & 0 & 1 & 0 & 1 & 3 & 1 & 1 & 0 & 1 & 1 & 1 \\
 0 & 1 & 0 & 0 & 0 & 1 & 0 & 0 & 0 & 1 & 0 & 0 & 0 & 0 & 0 \\
 0 & 1 & 0 & 0 & 0 & 0 & 0 & 1 & 0 & 0 & 1 & 0 & 0 & 0 & 0 \\
 1 & 0 & 2 & 0 & 1 & 0 & 1 & 0 & 0 & 0 & 0 & 1 & 0 & 0 & 0 \\
 0 & 0 & 0 & 1 & 0 & 0 & 0 & 1 & 0 & 0 & 0 & 0 & 1 & 0 & 0 \\
 0 & 0 & 0 & 1 & 0 & 1 & 0 & 0 & 0 & 0 & 0 & 0 & 0 & 1 & 0 \\
 0 & 1 & 0 & 1 & 0 & 1 & 0 & 1 & 2 & 1 & 1 & 0 & 1 & 1 & 1
 \end{array} \right]
 \end{equation*}
The initial vector is as follows:
\[ v = (1, 0, 0, 0, 0, 0, 0, 1, 0, 0, 0, 0, 0, 0, 0) \]
We then obtain
\begin{align*}
(Av)\cdot d_2 &= 12 \\
(A^2v)\cdot d_2 &= 75 \\
(A^3v) \cdot d_2 &= 384 \\
\vdots
\end{align*}
which yields the sequence of the number of spanning trees on $C_2(n)$.

Similar to the process with grids, there is nothing specific here to the simple cylinder graph - these methods can be used to obtain sequences for graph families of the form $G \times C_n$ for arbitrary $G$. However, due to the rapid growth of $B_{2k}$, the ability to find the appropriate matrices becomes somewhat impossible starting at graphs with five vertices. Nevertheless, we still have the following:

\begin{theorem}
For a given graph $G$ on $k$ vertices, there is a $B_{2k} \times B_{2k}$ matrix $M$ and a vector $v$ of length $B_{2k}$ such that
\[ (M^nv)\cdot d_k \]
is the number of spanning trees of the graph $G \times C_n$.
\end{theorem}

\begin{corollary}
For a given graph $G$ on $k$ vertices, the number of spanning trees $\set{T_n}$ of $G \times C_n$ satisfies a linear recurrence of order at most $B_{2k}$.
\end{corollary}

Although the sequence for $C_2(n)$ is already known, these methods used were able to obtain sequences for $C_3(n)$ and $K_3 \times C_n$, which we now state:
\begin{center}
\begin{tabular}{|l|}
\hline
$C_2(n):$ (\cite{FaaseSite}, with improvements) \\
$T_n = 10T_{n-1} - 35T_{n-2} + 52T_{n-3} - 35T_{n-4} + 10T_{n-5} - T_{n-6}$ \\
Sequence: $\set{1, 12, 75, 384, 1805, \ldots}$ (OEIS A006235)\\
Generating Function: $\frac{x (1 + 2 x - 10 x^2 + 2 x^3 + x^4)}{(-1 + 5 x - 5 x^2 + x^3)^2}$\\
\hline
$C_3(n):$ (new) \\
$T_n = 48T_{n-1} - 960T_{n-2} + 10622T_{n-3} - 73248T_{n-4} + 335952T_{n-5} - 1065855T_{n-6} + 2396928T_{n-7}$ \\
$-3877536T_{n-8} + 4548100T_{n-9} - 3877536T_{n-10} + 2396928T_{n-11} - 1065855T_{n-12} + 335952T_{n-13}$ \\
$-73248T_{n-14} + 10622T_{n-15} - 960T_{n-16} + 48T_{n-17} - T_{n-18}$ \\
Sequence: $\set{1, 70, 1728, 31500, 508805, \ldots}$ (OEIS to be submitted)\\
Generating Function: See \cite{RaffSite} \\
\hline
$K_3 \times P_n:$ (new) \\
$T_n = 58T_{n-1} - 1131T_{n-2} + 8700T_{n-3} - 29493T_{n-4} + 43734T_{n-5} - 29493T_{n-6} + 8700T_{n-7} $\\
$- 1131T_{n-8} + 58T_{n-9} - T_{n-10}$ \\
Sequence: $\set{3, 318, 12960, 410700, 11870715, \ldots}$ (OEIS to be submitted)\\
Generating Function: $\frac{3 x (1 + 48 x - 697 x^2 - 2474 x^3 + 9918 x^4 + 62 x^5 - 2045 x^6 +
   96 x^7 + 5 x^8)}{(-1 + 29 x - 145 x^2 + 145 x^3 - 29 x^4 + x^5)^2}$ \\
\hline
\end{tabular}
\end{center}

\section{Acknowledgements}

Special thanks to Andrew Baxter for thoroughly reviewing the paper and suggesting many helpful additions. Thanks also to Prof. Doron Zeilberger for taking me on as his seventh concurrent student, even though his self-proclaimed limit is four.

\bibliographystyle{plain}
\bibliography{refs}

\end{document}